\documentclass[11pt]{amsart}
\usepackage{amssymb,amsfonts,epsfig}
\title[Morse theory, Floer theory and closed geodesics of flat tori]
  {Morse theory on the loop space of flat tori and symplectic Floer theory}
\author{Joa Weber}
\address{Institute for Mathematical Sciences,
  SUNY at Stony Brook, NY 11794-3651, USA}
\email{joa@math.sunysb.edu, http://www.math.tu-berlin.de/$^\sim$joa}
\date{December 1996, revised March 1999}
\newtheorem{thm}{Theorem}[section]
\newtheorem{lem}[thm]{Lemma}

\theoremstyle{definition}
\newtheorem{defn}[thm]{Definition}
\newtheorem{remk}[thm]{Remark}

\newcommand{\1}{{{\mathchoice {\rm 1\mskip-4mu l} {\rm 1\mskip-4mu l}
  {\rm 1\mskip-4.5mu l} {\rm 1\mskip-5mu l}}}}
\newcommand{\cal}{\mathcal}

\newcommand{\s}{\scriptstyle}

\newcommand{\R}{{\mathbb R}}

\begin{document}                                      

\begin{abstract}
  We use closed geodesics to construct and compute 
  Bott-type Morse homology groups for the 
  energy functional on the loop space
  of flat $n$-dimensional tori, $n\ge 1$, and
  Bott-type Floer cohomology groups for their cotangent bundles 
  equipped with the natural
  symplectic structure. Both objects are isomorpic
  to the singular homology of the loop space.
  In an appendix we perturb the equations
  in order to eliminate degeneracies and to get to
  a situation with nondegenerate critical points only.
  The (co)homology groups turn out
  to be invariant under the perturbation.
\end{abstract}


\maketitle
\tableofcontents
 
  \section{Introduction}

  We are going to use the set of critical points of the symplectic
  action functional
  \begin{eqnarray} \label{17}
    {\cal A}_{H_{0}}\; : \; \Lambda T^{*}T^{n} & \rightarrow & {\mathbb R} \\
    z & \mapsto & \int_{{\mathbb R}/{\mathbb Z}} z^{*}\Theta - \int_{0}^{1}
    H_{0}(z {\scriptstyle (} t {\scriptstyle )} )\: dt \; , \nonumber
  \end{eqnarray}
  to construct an algebraic chain complex whose homology represents the
  singular homology of the free loop space of $T^n$.
  In (\ref{17}) $\Theta$ is the Liouville form on $T^*T^n$ and
  the Hamiltonian $H_0\: :\: T^*T^n \rightarrow {\mathbb R}$ is given by
  $z \mapsto \frac{1}{2} |z|^2$,
  where the metric $g$ on $T^n\subset{\mathbb C}^n$ is induced by the real part of the 
  hermitian inner product on ${\mathbb C}^n$.
  In our construction of the chain complex we follow the ideas in \cite{RT} and \cite{AB}.
  Because $H_0$ is time independent, the set $Crit\: {\cal A}_{H_0}$ of critical points of 
  ${\cal A}_{H_0}$ cannot be discrete: Assume a critical point $z$ is a nonconstant loop,
  then $S^1 \ni \tau \mapsto z(\cdot +\tau)$ produces a $S^1$-family of critical points.
  Therefore standard nondegenerate Floer theory does not apply. 
  We have to take into account the singular \cite{RT}
  or de Rham \cite{AB} (co)homology of $Crit\: {\cal A}_{H_0}$ leading
  to a double complex.

    In the general case the $L^2$-gradient flow lines of ${\cal A}_{H_0}$ between
  different connected components $Crit\: {\cal A}_{H_0}$ will
  enter the construction, but here they cannot exist
  as different components lie in different connected components of $\Lambda T^*T^n$. 
  Therefore the Bott-type Floer complex reduces to the singular chain
  complex of the critical manifolds, whose components are diffeomorphic to $T^n$,
  and so
  \begin{equation} \label{eq:Floer-Bott}
    HF^{-i}_{a,Bott}(T^{*}T^{n},H_0,g;{\mathbb Z}) =
    \begin{cases}
      \bigoplus_{l \in {\mathbb Z}^n , |l|\le |k|} {\mathbb Z}^{\binom{n}{i}}
       & , -i=-n,\ldots ,0 \\
      0 & ,\text{else}
    \end{cases}
  \end{equation}
  where $a=2\pi^2 |k|^2$ indicates that we consider only loops in $T^*T^n$
  whose symplectic action is less or equal to $a$. This set is denoted
  by $\Lambda^a T^*T^n$. The direct sum in
  (\ref{eq:Floer-Bott}) parametrizes the connected components of $\Lambda^a T^*T^n$.
  The grading $-i$ is provided by a generalized Conley-Zehnder index
  associated to any element of $Crit\: {\cal A}_{H_0}$ plus half the
  local dimension of $Crit\: {\cal A}_{H_0}$.

    On the other hand one can perturb $H_0$ by a time dependent function 
  (potential) $V$ on $T^n$, such that $Crit\: {\cal A}_{H_0+V}$ consists
  of isolated critical points and in this case standard Floer theory applies
  (this is the approach in \cite{SW} and
  \cite{W99}).
  The grading then is given by the standard Conley-Zehnder index \cite{CZ}
  and the trajectories of the negative gradient flow of ${\cal A}_{H_0+V}$
  connect critical points of \emph{increasing} nonpositive indices.
  That is why the cohomology notation is used.

    It is known that the Hamiltonian flow on $T^*T^n$ generated by 
  $H_0$ gives rise to geodesics
  of $(T^n,g)$ when projected to the base $T^n$. In this way
  $Crit\: {\cal A}_{H_0}$ corresponds precisely to the set of closed
  geodesics of $(T^n,g)$, which are the critical points of
  the classical action functional on the free loop space of $T^n$
  (called energy functional in Riemannian geometry)
  \begin{eqnarray} \label{18}
    {\mathcal I}_{0}\; : \; \Lambda T^{n} & \rightarrow & {\mathbb R} \\
    \gamma & \mapsto & {\scriptstyle \frac{1}{2} } \int_{0}^{1}
    |\dot \gamma (t)|^2 \: dt\; . \nonumber
  \end{eqnarray} 
  As the critical points of the symplectic and the classical action
  are naturally identified we denote them simply by $Crit$.
  Observe that on $Crit$ both functionals agree. This
  continues to hold in the presence of a time-dependent perturbation term.
  We use closed geodesics $\gamma$ with ${\mathcal I}_{0} (\gamma) \le a$ 
  to construct Bott-type Morse
  chain groups graded by the Morse index.
  The resulting homology groups $HM_{i}^{a,Bott}(
  \Lambda T^n,{\mathcal I}_{0},g;{\mathbb Z})$ turn out to be isomorphic 
  to $HF^{-i}_{a,Bott}(T^{*}T^{n},H_0,g;{\mathbb Z})$ as well as to
  the $i^{th}$ singular integral homology groups of 
  $$
    \Lambda^a T^n=\{ \gamma \in 
  \Lambda T^n \mid {\mathcal I}_{0} (\gamma)\le a \} ,
  $$
  which we compute using classical Morse theory \cite{Kl}.

  Recently  Viterbo \cite{Vi} as well as Salamon and the present author \cite{SW}, 
  \cite{W99} proved
  that singular homology of the loop space of a compact oriented Riemannian manifold
  is isomorphic to Floer cohomology of its cotangent bundle, hence the homology groups
  calculated in appendix \ref{se:perturbations}
  do in fact not depend on the metric or the potential.
  The former proof uses generating function homology, whereas
  in the second proof one defines integral Morse homology
  of the classical action functional, which is isomorphic to singular homology 
  of the free loop space,
  and then constructs a natural isomorphism to
  integral Floer homology of the cotangent bundle by showing a 1-1 correspondence
  of flow trajectories.

  One consequence of the perturbation of $H_0$ by the time dependent
  potential $V$ is that the critical
  points of ${\mathcal I}_V$ cannot be interpreted as closed geodesics any
  more. However, the corresponding homology groups should generally
  be the same. Following \cite{CFHW} we construct such a perturbation
  in the case of the $1$-sphere $S^1$ explicitly. As one expects we'll see that a
  connected component of $Crit\: {\cal A}_{H_0}$ resp.  $Crit\: {\mathcal I}_{0}$
  splits in a number of isolated critical
  points (depending on the perturbation) of Conley-Zehnder index $-1$ and $0$ resp.
  Morse index $1$ and $0$, such that the corresponding
  local homology groups are isomorphic to $H_{*}^{sing}(S^1;{\mathbb Z}_2)$.
  The case of more general critical manifolds as $S^1$ will be subject of
  future research. In what follows
  we denote $H_0$ by $H$, ${\mathcal I}_{0}$ by ${\mathcal I}$ 
  and ${\cal A}_H$ by ${\cal A}$.
  Moreover, we will use throughout Einstein's summation convention.

  {\bf Acknowledgements} 
  {\footnotesize I would like to thank Helmut Hofer for drawing
  my attention to this field. 
  For valuable discussions and support
  I thank my friend and colleague Kai Cieliebak as well as my advisor
  Ruedi Seiler. I am grateful to Graduiertenkolleg
  ``Geometrie und nichtlineare Analysis'' at Humboldt-Universit\"at
  zu Berlin for financial support.}

  \section{Bott-type Morse theory on the loop space} \label{15}

    Let $S^{1}$ be embedded in ${\mathbb C}$ as the unit circle
    and $i:T^n=S^{1}\times \ldots \times S^1 \hookrightarrow
    {\mathbb C}^n$. By $g$ we denote the flat riemannian metric on $T^n$ inherited from
    the real part of
    the hermitian inner product $\langle x,y \rangle_{{\mathbb C}^n}=\sum_{j=1}^n \bar{x}_j y_j$
    on ${\mathbb C}^n$; $\nabla$ denotes
    the Levi-Civita connection of $g$. With respect to the
    natural parametrization 
    $\psi^{-1} : {\mathbb R}^n / {\mathbb Z}^n \rightarrow T^n \subset {\mathbb{C}}^n$, 
    $(u^1,\ldots ,u^n) \mapsto (e^{2\pi iu^1} , \ldots ,e^{2\pi iu^n})$, 
    $g$ is given by $g=g_{jk} \: du^j\otimes du^k=(2\pi)^2 \delta_{jk} \: du^j\otimes du^k$ 
    and the volume element by $vol_{g}=(2\pi)^n \:
    du^1 \ldots du^n$.

    The free loop space $\Lambda T^n=H^1 ({\mathbb R} / {\mathbb Z},T^n)$ is defined
    to be the completion of $C^\infty ({\mathbb R} / {\mathbb Z},T^n)$ with respect
    to the norm on $C^\infty ({\mathbb R} / {\mathbb Z},{\mathbb C}^n)$
    $$
      \| \gamma \|_{1,2}^2 =\int_0^1 \langle \gamma (t),
      \gamma (t) \rangle_{{\mathbb C}^n} dt 
      + \int_0^1 \langle \partial_t \gamma (t), 
      \partial_t \gamma (t) \rangle_{{\mathbb C}^n} dt .
    $$
    Note that for $\gamma=(e^{2\pi i \gamma^1} , \ldots ,e^{2\pi i \gamma^n}),
    \tilde{\gamma} \in C^\infty ({\mathbb R} / {\mathbb Z},T^n)$ it turns out
    \begin{equation*}
      \begin{split}
      \| \gamma - \tilde{\gamma} \|_{1,2}^2
     &=2\int_0^1 \Bigl( n - \sum_{j=1}^n \cos 2\pi (\gamma^j - \tilde{\gamma}^j) \Bigr) dt \\
     &\quad +(2\pi)^2 \int_0^1 \sum_{j=1}^n \Bigl( (\dot \gamma^j)^2 
      + (\dot{\tilde{\gamma}}^j)^2
      -2 \dot \gamma^j \dot{\tilde{\gamma}}^j \cos 2\pi (\gamma^j - \tilde{\gamma}^j) \Bigr) 
      dt .
      \end{split}
    \end{equation*}
    By the Sobolev embedding theorem $\Lambda T^n$ embedds into $C^0 ({\mathbb R}/{\mathbb Z},T^n)$.
    In what follows we occasionally identify $T^n$ with ${\mathbb R}^n/{\mathbb Z}^n$.

    \begin{lem}
      (cf. \cite{Jo} Se. 5.4) The energy functional ${\mathcal I}$ 
      is continuously differentiable.
    \end{lem}

    \begin{lem}
      (cf. \cite{Jo} Lemma 7.2.1 \& Se. 5.4, \cite{Kl} Thm. 1.3.11)
      The critical points of ${\mathcal I}$ in $\Lambda T^n$ are precisely the
      closed geodesics of $(T^n,g)$.
    \end{lem}

    The Christoffel symbols of $\nabla$ vanish
    because the matrix elements $g_{jk}$ are constant and so
    $(\gamma^1 ,\ldots, \gamma^n) \in C^{\infty}({\mathbb R}/{\mathbb Z},{\mathbb R}^n /{\mathbb Z}^n )$ 
    is a closed geodesic, if and only if
    for $j=1,\ldots,n$
    \begin{equation} \label{1}
      {\textstyle \frac{d^{2}}{dt^{2}}} \gamma^j (t) = 0 .
    \end{equation}
    Now $\gamma^j$ considered as element of 
    $C^{\infty}({\mathbb R} ,{\mathbb R} )$ is a solution of
    (\ref{1}) if and only if
    \begin{displaymath}
      \gamma^j (t)=v^j t+q^j .
    \end{displaymath}
    The condition $\gamma^j (t+1)=\gamma^j (t) \pmod{1}$
    implies $v^j=k^j\in {\mathbb Z} $, therefore
    \begin{displaymath}
      Crit\: {\mathcal I}=\{ \gamma_{k,q} (t)=k t+q 
      \mid k \in {\mathbb Z}^n,
      q \in {\mathbb R}^n / {\mathbb Z}^n \} .
    \end{displaymath}
    There is a grading on $Crit\: {\mathcal I}$ given by the Morse index.
    \begin{lem}
      (cf. \cite{Jo} Thm. 4.1.1 \& Se. 5.4)
      Let $\gamma \in Crit\: {\mathcal I}$ and 
      $X,Y\in C^\infty({\mathbb R}/{\mathbb Z},{\mathbb R}^n)$, i.e. smooth
      vector fields along $\gamma$, then
      \begin{displaymath}
        d^{2}{\mathcal I} (\gamma )\ (X,Y)=\int_{0}^{1} 
        g(\dot X(t),\dot Y(t))\ dt\ .
      \end{displaymath}
    \end{lem}
    The associated selfadjoint operator is the \emph{Jacobi operator}      
    \begin{eqnarray*}
      L=L_\gamma : L^2({\mathbb R}/{\mathbb Z},{\mathbb R}^n) \supset 
      H^2({\mathbb R}/{\mathbb Z},{\mathbb R}^n) & \to & L^2({\mathbb R}/{\mathbb Z},{\mathbb R}^n) \\
      X & \mapsto & -{\textstyle \frac{d^{2}}{dt^{2}}} X  .
    \end{eqnarray*}
    For $j=1,\ldots,n$ the solution of
    \begin{equation} \label{2}
      -{\textstyle \frac{d^{2}}{dt^{2}}} X^j(t)=0
    \end{equation}
    is, as an element of $C^{\infty}({\mathbb R} ,{\mathbb R} )$, given by $X^j(t)=v^jt+\xi^j$, 
    with $v^j, \xi^j \in {\mathbb R} $,
    but now the periodicity condition $X^j(t+1)=X^j(t)$ implies $v^j=0$ and so
    \[
      ker \: L\simeq \{ X (t)=\xi \mid \xi \in {\mathbb R}^n \} \simeq {\mathbb R}^n  .
    \]

    \begin{remk}
      Note that for any Riemannian manifold $(M,g)$ if $\gamma$ is a closed
      geodesic, then $\dot \gamma \in Ker\: L$. Therefore
      the kernel of the Hessian $d^{2}{\mathcal I} (\gamma )$ is always at 
      least 1-dimensional.
      In order to get a trivial kernel one has to introduce a time-dependent
      perturbation of ${\mathcal I}$.
      This will be carried out in section \ref{se:perturbations}.
    \end{remk}

    Next we are interested in the negative eigenvalues of $L$. These do not
    exist, because the operator $-d^{2}/dt^{2}$ is positive semidefinite
    on $C^{\infty}({\mathbb R}/{\mathbb Z},{\mathbb R} )$. The reason is the periodicity of
    the domain ${\mathbb R}/{\mathbb Z}$. To see this Fourier decompose $X (t)$ and
    apply $-d^{2}/dt^{2}$ to each summand. More generally, this follows by partial
    integration and the closedness of the manifold.
    
    \begin{defn}
      For $\gamma \in Crit\: {\mathcal I}$ we define its \emph{Morse index}
      and \emph{nullity} to be the number of negative eigenvalues
      of $L$ (counted with multiplicities) and the
      dimension of its kernel, respectively.
    \end{defn}
    So in our case we have $Ind_{{\mathcal I}}(\gamma )=0$ and $Null_{{\mathcal I}}(\gamma )=n$
    for all $\gamma \in Crit\: {\mathcal I}$ and the results derived so far may be
    summarized as follows.
    \begin{lem}
      For $(T^n,g=\iota^* \langle \cdot , \cdot \rangle_{{\mathbb C}^n})$ as above we have
      \begin{equation*}
      \begin{gathered}
        Crit \:{\mathcal I} =\bigsqcup_{k\in {\mathbb Z}^n } G^{k}\; , \\
        G^{k}=\{ \gamma_{k,q}(t)=(e^{2\pi i(k^1t+q^1)},\ldots, e^{2\pi i(k^nt+q^n)})
        \mid q \in {\mathbb R}^n /{\mathbb Z}^n \} \; ,
      \end{gathered}
      \end{equation*}
      i.e. $G^{k}$ is a submanifold of $\Lambda T^n$ diffeomorphic to $T^n$.
      $G^{k}$ is a \emph{nondegenerate critical submanifold in the sense of Bott}, 
      i.e. $d^{2}{\mathcal I}(\gamma )$
      restricted to the normal bundle $N_{\gamma}G^{k}$ is nondegenerate
      for any $\gamma \in G^{k}$.
    \end{lem}

    The last statement follows, because one can canonically identify
    $ker \: L_\gamma \simeq ker \: d^2 {\mathcal I} (\gamma)$
    with $T_\gamma G^k$ by 
    $$
      X (t)=\xi \mapsto \xi ={\textstyle \left. \frac{d}{dt} \right|_{t=0} }
      (kt+q+\tau \xi) .
    $$
 
    The critical submanifolds generate an algebraic chain  complex
    as follows: 
    The $i^{th}$ \emph{Bott-type Morse chain group}
    is defined to be the $i^{th}$ singular chain group of 
    $Crit^a =Crit\: {\mathcal I} \cap \Lambda^aT^n$, $a=2\pi^2 |k|^2$,
    with coefficients in ${\mathbb Z}$

    \begin{displaymath}
      CM_{i}^{a,Bott}(\Lambda T^n ,{\mathcal I},g;{\mathbb Z})=
      \begin{cases}                                                   
        C^{sing}_{i}\left( \bigsqcup_{l\in {\mathbb Z}^n,|l|\le|k|} G^l ; {\mathbb Z} \right) 
        & ,\: i\in {\mathbb N}_{0} \: ,\\
        0 & ,\: i\in {\mathbb Z}\setminus {\mathbb N}_{0} \: .              
      \end{cases}
    \end{displaymath}   

    If $k\not=j$, then $G^k$ and $G^j$ lie in different
    connected components of $\Lambda T^n$,
    so we cannot expect to
    have any connecting orbit (in the sense of Morse/Floer theory)
    between $G^k$ and $G^j$.
    Therefore we may define the
    \emph{Bott-type Morse boundary operator} $d_{i}^{M}$ simply to be the singular
    boundary operator $\partial_{i}^{sing}$ on the singular chains of our critical manifolds
    \begin{eqnarray*}
      d_{i}^{M} : CM_{i}^{a,Bott}(\Lambda T^n,{\mathcal I},g;{\mathbb Z}) & \to &
      CM_{i-1}^{a,Bott}(\Lambda T^n,{\mathcal I},g;{\mathbb Z}) \; , \\ 
      x & \mapsto & \partial_{i}^{sing} x\: .
    \end{eqnarray*}
    The chain complex property $d^{M}_{i-1}\circ d^{M}_{i}=0$ for all $i\in {\mathbb Z}$
    then follows trivially from the one of the singular chain complex and so
    we may define the \emph{Bott-type Morse homology groups} to be
    \begin{displaymath}
      HM_{i}^{a,Bott}(\Lambda T^n,{\mathcal I},g;{\mathbb Z})
      =\frac{ker\; d^{M}_{i}}{im\; d^{M}_{i+1}},
      \; i\in {\mathbb Z} \; .
    \end{displaymath}
    For $i =0,\ldots,n$ a simple computation gives (the groups are $0$ else)
    \begin{eqnarray} \label{eq:Morse-hom}
      &&HM_{i}^{a,Bott}(\Lambda T^n,{\mathcal I},g;{\mathbb Z})  = H_{i}^{sing}
       (\bigsqcup\nolimits_{l\in {\mathbb Z}^n,|l|\le |k|} G^l ; {\mathbb Z}) \\
      && \; \; \; \; \; \; \; \; \; = H_{i}^{sing}
       (\bigsqcup\nolimits_{l\in {\mathbb Z}^n,|l| \le |k|} T^n ; {\mathbb Z})  
       = \bigoplus_{l\in {\mathbb Z}^n,|l| \le |k|} {\mathbb Z}^{\binom{n}{i}} .
    \end{eqnarray}

  \section{Singular homology of the loop space}

    By \cite{Kl} Thm. 1.2.10 the inclusion $\Lambda T^n=H^1
    ({\mathbb R}/{\mathbb Z},T^n)\hookrightarrow C^{0}({\mathbb R}/{\mathbb Z}
    ,T^n)$ is a homotopy equivalence, hence
    \begin{displaymath}
      H_{*}^{sing}(\Lambda T^n;{\mathbb Z})\simeq
      H_{*}^{sing}(C^{0}({\mathbb R}/{\mathbb Z},T^n);{\mathbb Z})\; .
    \end{displaymath}
    The set $\pi_{1}(T^n)$ of free homotopy classes of continuous maps
    equals ${\mathbb Z}^n$. A homotopy $F_{s}$, $s\in [0,1]$, between 
    $\gamma , \tilde{\gamma} \in C^{0}({\mathbb R}/{\mathbb Z},T^n)$, may be considered as a path in
    $C^{0}({\mathbb R}/{\mathbb Z},T^n)$ from $\gamma$ to $\tilde{\gamma}$. Hence the homotopy classes
    $\pi_1 (T^n)$ correspond precisely to the pathwise connected components of
    $C^{0}({\mathbb R}/{\mathbb Z},T^n)$. On the other hand to any pathwise connected
    component of $C^{0}({\mathbb R}/{\mathbb Z},T^n)$ corresponds a
    generator of $H_0^{sing}(C^{0}({\mathbb R}/{\mathbb Z},T^n);{\mathbb Z})$,
    therefore
    \begin{displaymath}
      H_0^{sing}(\Lambda T^n;{\mathbb Z})\simeq
      H_0^{sing}(C^{0}({\mathbb R}/{\mathbb Z},T^n);{\mathbb Z})=
      \bigoplus_{\alpha \in \pi_{1}(T^n)} {\mathbb Z}
      =\bigoplus_{k\in {\mathbb Z}^n} {\mathbb Z} \; .
    \end{displaymath}
    We are going to compute the higher homology groups via classical Morse theory
    of the energy functional ${\mathcal I}$ on the loop space $\Lambda T^n$.
    As we have seen in the former section the critical submanifolds $G^k$ of 
    $\Lambda T^n$ with respect to ${\mathcal I}$ are nondegenerate
    in the sense of Bott, they have Morse index $0$ and they are diffeomorphic
    to $T^n$.
    $G^0$ corresponds to the trivial (constant) geodesics and
    for any element of $G^k$
    \[
      \gamma_{k,q} =( e^{2\pi i (k^1t+q^1)} ,\ldots, e^{2\pi i (k^nt+q^n)})
    \]
    we compute
    \begin{displaymath}
      {\mathcal I} (\gamma_{k,q}) = {\textstyle \frac{1}{2}}
      \int_0^1 \langle \partial_t \gamma_{k,q} , \partial_t \gamma_{k,q} \rangle_{{\mathbb C}^n} dt
      =2 \pi^2 |k|^2 .
    \end{displaymath}

    \begin{thm}
      (cf. \cite{Kl}, Thm. 2.4.10) Let $(M,g)$ be a compact Riemannian
      manifold and assume that the set of critical points of ${\mathcal I}$
      in ${\mathcal I}^{-1}(a)$ is a nondegenerate critical submanifold $B$
      of $\Lambda M$. Then there exists $\epsilon >0$, such that
      $\Lambda^{a+\epsilon}M$ is (equivariantly) diffeomorphic
      to $\Lambda^{a-\epsilon}M$ with closed disk bundle of type
      $D\mu^{-} \oplus D\mu^{+}$ attached. $\mu$ denotes the normal
      bundle of $B$ and $\mu =\mu^{-} \oplus \mu^{+}$ is the 
      decomposition in the negative and positive subbundle (w.r.t.
      the Hessian of ${\mathcal I}$, i.e. $rk\: \mu^- =Ind_{\mathcal I} (\gamma)$
      for $\gamma \in B$).
    \end{thm}

    In our case all critical submanifolds have Morse index $0$ and this
    implies that for $\epsilon >0$ sufficiently small
    \begin{displaymath}
      \Lambda^{2\pi^{2} |k|^{2}+\epsilon } T^n \simeq
      \bigsqcup_{l\in {\mathbb Z}^n , |l|\le |k|} D\mu_{l}^{+}\; .
    \end{displaymath}
    Any of the bundles $D\mu_{l}^{+} \rightarrow G^l$ is contractible
    on $G^l$, hence
    \begin{equation*}
    \begin{split}
       HM_{i}^{sing}(\Lambda^{2\pi^{2} |k|^{2}+\epsilon } T^n;{\mathbb Z})
      &\simeq H_{i}^{sing}
       (\bigsqcup\nolimits_{l\in {\mathbb Z}^n , |l|\le |k|} T^n ; {\mathbb Z}) \\
      &= 
       \begin{cases}                                 
         \bigoplus_{l\in {\mathbb Z}^n , |l|\le |k|} 
         {\mathbb Z}^{\binom{n}{i}} & , i=0, \ldots , n \\
         0 & , \text{else}\: .                     
       \end{cases}
    \end{split}
    \end{equation*}

  \section{Bott-type Floer homology of the cotangent bundle}

    Let $(T^n\subset {\mathbb C}^n , g)$
  be as above, then the parametrization $\psi : {\mathbb C}^n \supset
  T^n \to {\mathbb R}^n / {\mathbb Z}^n$, $(e^{2\pi i u^1},\ldots,e^{2\pi i u^n}) 
  \mapsto (u^1,\ldots,u^n)$
  induces natural coordinates on $TT^n$ and $T^*T^n$
  \begin{equation*}
    \begin{split}
      d\psi : {\mathbb C}^n \times {\mathbb C}^n \supset TT^n 
     &\to ({\mathbb R}^n / {\mathbb Z}^n) \times {\mathbb R}^n \\
      (e^{2\pi i u^1},\ldots,e^{2\pi i u^n},
      2\pi iv^1 e^{2\pi i u^1},\ldots,2\pi iv^n e^{2\pi i u^n})
     &\mapsto (u^1,\ldots,u^n,v^1,\ldots,v^n)
    \end{split}
  \end{equation*}
  and
  \begin{equation*}
    \begin{split}
      {d\psi^*}^{-1} : {\mathbb C}^n \times {{\mathbb C}^n}^* \supset T^*T^n 
     &\to ({\mathbb R}^n / {\mathbb Z}^n) \times {{\mathbb R}^n}^* \\
      (e^{2\pi i u^1},\ldots,
      {\textstyle \frac{-i}{2\pi}} v_1 e^{-2\pi i u^1},\ldots)
     &\mapsto (u^1,\ldots,u^n,v_1,\ldots,v_n)=(u,v) .
    \end{split}
  \end{equation*}
  Note that $\partial / \partial u^j$ is identified with $2\pi i e^{2\pi i u^j} e^j$
  and $du^j$ with $\frac{-i}{2\pi} e^{-2\pi i u^j} e^j$, where $e_j,e^j$ are
  elements of the standard bases of  ${\mathbb C}^n$ and ${{\mathbb C}^n}^*$ respectively.

    The free loop space $\Lambda T^* T^n=H^1({\mathbb R} / {\mathbb Z},T^*T^n)$
  is defined to be the completion of $C^\infty ({\mathbb R} / {\mathbb Z},T^*T^n)$
  with respect to the norm on $C^\infty ({\mathbb R} / {\mathbb Z},{\mathbb C}^n \times
  {{\mathbb C}^n}^*)$
  \begin{equation*}
    \begin{split}
      \| z \|_{1,2}^2 
     &=\| (x,y) \|_{1,2}^2 =\| (x,y) \|_2^2 + \| \partial_t (x,y) \|_2^2 \\
     &=\int_0^1 \Bigl( \langle x(t) ,x(t) \rangle_{{\mathbb C}^n}
      +\langle y(t) ,y(t) \rangle_{{{\mathbb C}^n}^*} \Bigr) dt \\
     &\quad +\int_0^1 \Bigl( \langle \dot x(t) , \dot x(t) \rangle_{{\mathbb C}^n}
      +\langle \dot y(t) , \dot y(t) \rangle_{{{\mathbb C}^n}^*} \Bigr) dt .
    \end{split}
  \end{equation*}

  Denote by $\tau_M :TM\rightarrow M$ the tangent -- and by $\tau_M^* :T^*M \rightarrow
  M$ the cotangent bundle of a manifold $M$.
  The \emph{Liouville form} $\Theta :TT^{*}T^n\rightarrow {\mathbb R}$ is defined by
  \begin{displaymath}
    \Theta (\zeta ) = (\tau_{T^{*}T^n} \zeta ) \: (T\tau^{*}_{T^n} \zeta ) ,
  \end{displaymath}
  i.e. in our local coordinates $(u,v)$ we have $\Theta \mid_{(u,v)} =v_j\:du^j$.
  The \emph{canonical symplectic form} on $T^{*}T^n$ is $\Omega=-d\Theta$,
  i.e. 
  \begin{displaymath}
    \Omega \mid_{(u,v)} (\cdot ,\cdot )=(du^j\wedge dv_j)\: (\cdot ,\cdot ) .
  \end{displaymath}
  Therefore $\Omega$ is represented by the standard symplectic form on ${\mathbb R}^{2n}$
  \begin{eqnarray} \label{11}
    \Omega_{0} \: : {\mathbb R}^{2n} \times {\mathbb R}^{2n} & \rightarrow & {\mathbb R} \\
    (x,\tilde{x} ) & \mapsto & (J_{0}x)^{T}\tilde{x} ,\nonumber
  \end{eqnarray}
  where
  \begin{displaymath}
    J_{0}=\begin{pmatrix} 0&-\1 \\ \1&0 \end{pmatrix}
  \end{displaymath}
  is the standard complex structure on ${\mathbb R}^{2n}$.
  
  To our Hamiltonian $H$ we assign
  the \emph{Hamiltonian vector field $X_{H}$} by setting
  \begin{equation} \label{12}
    dH(\cdot )=\Omega (X_{H},\cdot ) .
  \end{equation}

  We are interested in the critical points of the action functional ${\cal A}$
  because they are related to ${\cal P}er \: H$, 
  the time-$1$-periodic integral curves of the
  Hamiltonian vector field $X_{H}$, as follows
  \begin{displaymath}
    Crit\: {\cal A}=\{ x\in C^{\infty}({\mathbb R}/{\mathbb Z},T^{*}T^n)\: |\:
    \dot x=X_{H}(x)\} ={\cal P}er \: H \: .
  \end{displaymath}

  The differential of $H(u,v)=\frac{1}{2} |v|_g^2=
  \frac{1}{2} v_j v_k \delta^{jk} /(2\pi)^2$ is given by
  \begin{displaymath}
    dH(u,v)={\textstyle \frac{1}{(2\pi )^{2}}} \: v_j\: dv_j
    ={\textstyle \frac{1}{(2\pi )^{2}}}
    \begin{pmatrix} 0\\v \end{pmatrix} \; .
  \end{displaymath}

  The Hamiltonian vector field $X_{H}\: :\: T^{*}T^n\to TT^{*}T^n$
  is computed via its defining equation (\ref{12}): The lhs has been just
  calculated and the rhs will be determined by using the {\em Ansatz}
  $X_{H}(u,v)=(u,v;r^j {\scriptstyle (}u,v{\scriptstyle )} 
  \partial_{u^j} +s^k {\scriptstyle (}u,v{\scriptstyle )} \partial_{v_k} )\;$.
  By (\ref{11}) the quantity $\Omega_{(u,v)} (X_{H} {\scriptstyle (} u,v
  {\scriptstyle )} ,\cdot )$ is represented by
  \begin{displaymath}
    (J_{0}X_{H}{\scriptstyle (} u,v{\scriptstyle )})^{T} = (-s
    {\scriptstyle (}u,v{\scriptstyle )},r{\scriptstyle (}u,v{\scriptstyle )})\: .
  \end{displaymath}
  Comparing lhs and rhs of (\ref{12}) gives
  \[
    s(u,v)=0\; ,\; r(u,v)=v/(2\pi )^{2} \; ,\; \text{hence} \;
    X_{H}(u,v)={\textstyle \frac{1}{(2\pi )^{2}}} \begin{pmatrix} v \\ 0 \end{pmatrix} \; .
  \]

  The time-$1$-periodic trajectories of the Hamiltonian vector field
  $X_{H}$ are exactly the solutions to the initial value problem
   \begin{equation} \label{6}
    \begin{pmatrix} \dot u(t)\\ \dot v(t) \end{pmatrix} =
    {\textstyle \frac{1}{(2\pi )^{2}}} \begin{pmatrix} v(t)\\ 0 \end{pmatrix}\; \;  , \; \;
    \begin{pmatrix} u(0)\\ v(0) \end{pmatrix} =
    \begin{pmatrix} u_0\\ v^0 \end{pmatrix}  .
  \end{equation}
  The Ansatz
  \[
    x(t)=\begin{pmatrix} u(t)\\ v(t) \end{pmatrix} =
    \begin{pmatrix} \frac{v_0}{(2\pi )^{2}} t+u_0 \\ v^0 \end{pmatrix}
  \]
  solves (\ref{6}) for $x\in C^{\infty}({\mathbb R},{\mathbb R}^{2n})$.
  The condition $u\in C^{\infty}({\mathbb R}/{\mathbb Z},{\mathbb R}^n/{\mathbb Z}^n)$
  , i.e. $u(t+1)=u(t)\pmod{1}$, implies $v^0/(2\pi )^{2}=k\in {\mathbb Z}^n$, hence
  \begin{eqnarray} \label{7}
    {\cal C} &=& \{ x\in C^{\infty} ({\mathbb R}/{\mathbb Z},T^{*}T^n)|\dot x=X_{H}(x) \} \\
    {} & \simeq & \{ x_{u_0,k}(t)=
    \begin{pmatrix} {\scriptstyle kt+u_0} \\ {\scriptstyle (2\pi )^{2}k} \end{pmatrix}
    |\: k\in {\mathbb Z}^n ,u_0 \in {\mathbb R}^n /{\mathbb Z}^n \} \; .\nonumber
  \end{eqnarray}

  Now in Floer theory one assigns
  an integer $\mu_{CZ}$, called \emph{Conley-Zehnder index}
  (cf. \cite{CZ},\cite{SZ}), to any element $x=x_{u^0,k}\in {\cal C}$.
  Linearizing $\phi_{t}:{\mathbb R}^n/{\mathbb Z}^n \times {\mathbb R}^n \rightarrow 
  {\mathbb R}^n/{\mathbb Z}^n \times {\mathbb R}^n$, the
  time-$t$-map associated to $X_{H}$, leads to a path $A$ in $Sp(2n,\R)$
  \[
    \phi_{t} \begin{pmatrix} u_0 \\ v^0 \end{pmatrix} =x_{u_0,k}(t)=
    \begin{pmatrix} kt+u_0 \\ (2\pi )^{2}k \end{pmatrix}\; , \;
    A(t)=d\phi_{t} |_{(u_0,v^0)}=
    \begin{pmatrix} \1 & \frac{t}{(2\pi )^{2}} \1 \\ 0 & \1 \end{pmatrix}\; .
  \]
  As $A(0)=\1$, but $A(1)\notin Sp^{*}=\{ M\in Sp(2n,{\mathbb R}) \mid 1 \notin 
  spec(M)\}$, the usual definition of the Conley-Zehnder index does not apply.
  On the other hand for paths $\Psi$ starting at $\1$ and ending outside
  the Maslov cycle ${\cal C_+}=Sp(2n,\R) \setminus Sp^*$ it is shown in
  \cite{RS} remark 5.4 that
  \begin{equation} \label{eq:CZ}
    \mu_{CZ} (\Psi)=\mu_{Lag} (Graph \: \Psi,\Delta)
  \end{equation}
  where $\Delta \subset \R^{2n} \times \R^{2n}$ denotes the diagonal
  and $\mu_{Lag}$ is the Maslov index for \emph{any} continuous path of Lagrangian
  planes in $(\R^{2n} \times \R^{2n},-\omega_0 \oplus \omega_0)$.

  \begin{figure}[ht] \label{fig:A(t)}
    \centering\epsfig{figure=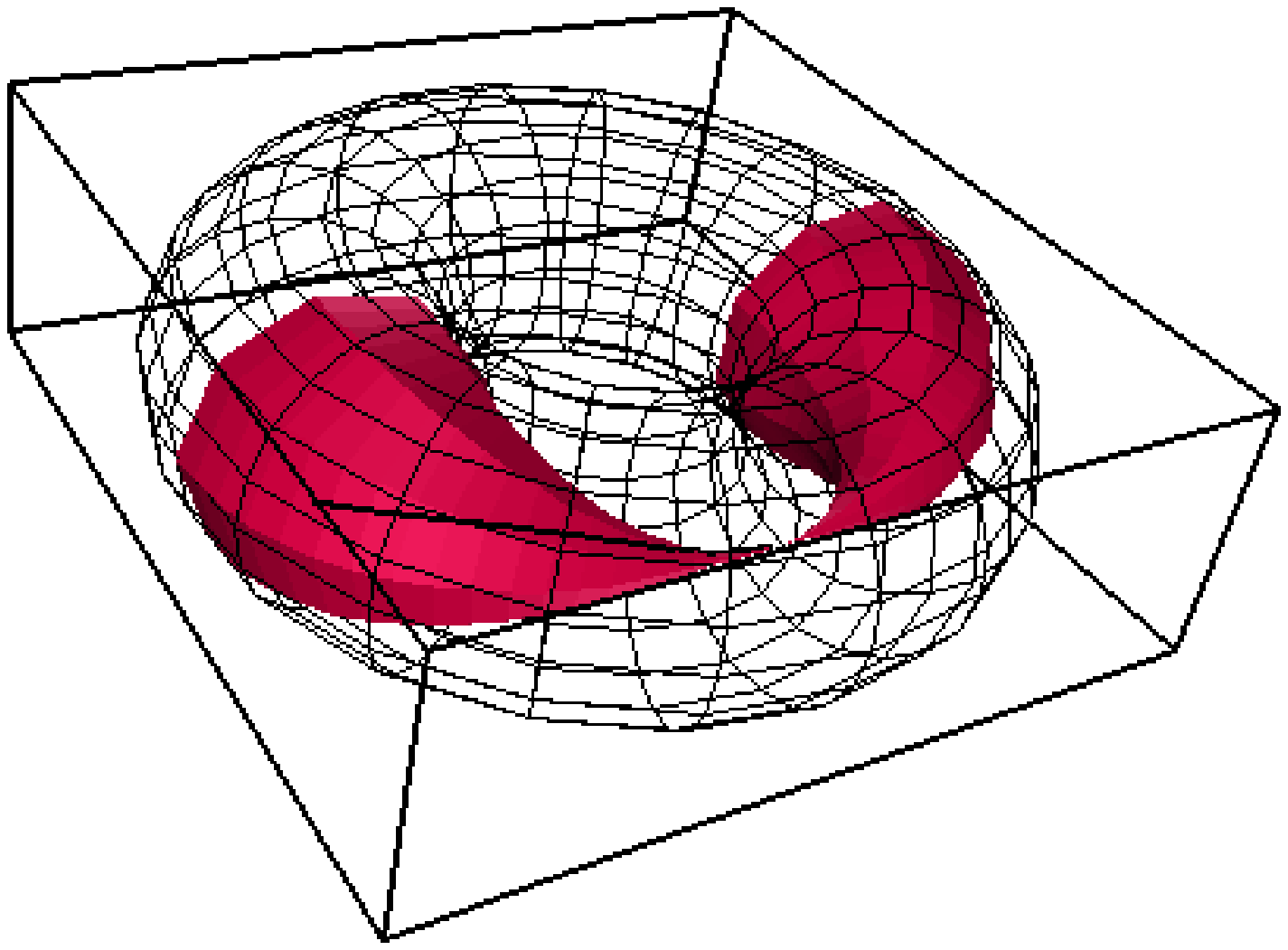,width=8cm} 
    \caption{The Maslov cycle $\cal C_+$ inside $Sp(2,\R)$}
  \end{figure}

    This index is invariant under homotopies of paths as long as the endpoints
  do not leave their stratum. It is therefore natural to define
  the \emph{generalized Conley-Zehnder index} for any continuous path $\Psi$ in
  $Sp(2n,\R)$ by the right hand side of (\ref{eq:CZ})
  \begin{equation} \label{eq:gen-CZ}
    \mu_{CZ}^g (\Psi)=\mu_{Lag} (Graph \: \Psi,\Delta) .
  \end{equation} 
  To show that $\mu_{CZ}^g (A)=-n/2$ it suffices -- in view of the
  product property of $\mu_{lag}$ -- to consider the case
  \[
    a(t)= \begin{pmatrix} 1&\frac{t}{4\pi^2} \\ 0&1 \end{pmatrix} .
  \]
  This path is shown in figure \ref{fig:A(t)} for $t \in [0,200]$
  together with the Maslov cycle (the codimension $1$ algebraic variety
  with one singular point, which corresponds to $\1$). Here $Sp(2,\R)$
  is identified with the open full torus using an explicit
  homeomorphism from \cite{GL}. Full details along with more pictures
  may be found in \cite{W98}, \cite{W99}.

    To calculate $\mu_{CZ}^g(a)$ we use its homotopy invariance.
  The polar decomposition $a(t)=S(t)R(t)$ gives a unitary path $R(t)=e^{2\pi it}$
  homotopic to $a$, where $a(0)=\1$ is preserved. Then connect $R(1)$ to
  $a(1)$ by a path $B(s)$ without intersecting $\cal C_+$ except for
  the endpoint $B(1)=a(1)$. Following first $R$ and then $B$ represents
  a path homotopic to $a$ with fixed endpoints. Near $a(1)$ the Maslov
  cycle is an embedding and we may use the intersection number interpretation
  of $\mu_{CZ}$ to get a contribution $+1$. At the singular point $a(0)=\1$ of
  $\cal C_+$ we compute the signature of the corresponding crossing form of the path $R$,
  cf. \cite{RS}, which turns out to be $-2$. As both endpoints lie in the Maslov
  cycle they are weighted by a factor $1/2$ and so
  \[
    \mu_{CZ}^g(a)=\frac{1}{2} (-2+1)=-\frac{1}{2} .
  \]
  As it should be, perturbing $a$ to a path $\tilde{a}$ on the other side
  of $\cal C_+$ leads to the same number
  \[
    \mu_{CZ}^g(\tilde{a})=\frac{1}{2} (0-1)=-\frac{1}{2} .
  \]

    So all elements of $Crit \: {\cal A}$ do have the same
  generalized Conley-Zehnder index $-n/2$ (note that this implies
  that there are no nonconstant trajectories of the negative gradient
  flow of $\cal A$ between critical points). Moreover, addition of 
  $\frac{1}{2} dim \: Crit \: {\cal A}=\frac{n}{2}$ to $\mu_{CZ}^g$
  yields the Morse index $0$ of the underlying closed geodesics.
  Let the \emph{Bott-type Floer cochain groups}
  be given by the singular chains in $Crit \: {\cal A}$, where the grading
  is minus the singular grading
  \[
    CF^{-i}_{a,Bott}(T^*T^n,H,g;{\mathbb Z} ) \; =\;
    \begin{cases}
      C_i^{sing}(Crit_a \: {\cal A}) & ,-i=-n,\ldots,0 \\
      0 & ,\mbox{else}
    \end{cases}
  \]
  where $a=2\pi^2 |k|^2={\cal A}(x_{u_0,k})$. This choice of
  grading is motivated by the nondegenerate case, 
  cf. appendix \ref{se:perturbations} and \cite{W99}, \cite{SW}.
  The coboundary operator $\delta_F^{-i}$ is defined to be
  the singular boundary operator $\partial_i^{sing}$ and so
  the \emph{Bott-type Floer cohomology groups}
  \begin{eqnarray*}
    HF^{-i}_{a,Bott}(T^{*}T^n,H,g;{\mathbb Z} ) &=& 
     \frac{ker\: \delta^{-i}_F}{im\: \delta^{-i-1}_F} \; ,\; i\in {\mathbb Z} \; .
  \end{eqnarray*}
  are given by
  \begin{eqnarray*}
    && HF^{-i}_{a,Bott}(T^*T^n,H,g;{\mathbb Z}) = H_{i}^{sing}(Crit_a \: {\cal A}) \\
    && \qquad \qquad=  H_{i}^{sing}
     \left(\bigsqcup\nolimits_{l\in {\mathbb Z}^n,|l| \le |k| } T^n ; 
     {\mathbb Z}_{2}\right) \\
    && \qquad \qquad=
     \begin{cases}                                 
       \bigoplus_{l\in {\mathbb Z}^n,|l| \le |k|} 
       {\mathbb Z}^{\binom{n}{i}}  & ,\; -i=-n,\ldots,0 \\
       0 & ,\; \text{else}\: .                     
     \end{cases}
  \end{eqnarray*}

\begin{appendix}
\section{Perturbations} \label{se:perturbations}

  For simplicity let us consider the case $n=1$ only. The general case
then follows by taking product manifolds and direct sums of Morse functions.
Moreover, we restrict to the $k^{th}$ connected component $\Lambda_kS^1$ resp. 
$\Lambda_kT^*S^1$ of the loop space consisting of loops of winding number
$k$. This restriction will be clear from our method of perturbation,
which does not work uniformly for all components.

First we are going to destroy the circle degeneracy of a closed geodesic
$\gamma_{k,q_0} (t)=kt+q_0$ (represented as a map from ${\mathbb R}/{\mathbb Z}$
to itself). We follow the beautiful construction in \cite{CFHW}. The strategy
is as follows: Choose a Morse-function $h$ on $\gamma_{k,q_0} ({\mathbb R}/{\mathbb Z})
={\mathbb R}/{\mathbb Z}$, e.g. $h(y)=-\cos 2\pi y$, define a time-dependent potential
$V$ on $\gamma_{k,q_0} ({\mathbb R}/{\mathbb Z})$
\[
  V(t,\gamma_{k,q_0} {\s (}s{\s )}):=h(ks-kt)=-\cos 2\pi (\gamma_{k,q_0} {\s (}s{\s )}
  -kt-q_0)\; .
\]
The critical points of the perturbed classical action functional ${\cal I}_V(\gamma )=
\int_0^1 \frac{1}{2} \mid \dot \gamma \mid^2 -V(t,\gamma ) \: dt$ are
solutions of
\begin{eqnarray}\label{eq:perM}
  0 &=& - \stackrel{\cdot \cdot}{\gamma} (t)-\nabla V(t,\gamma {\s (}t{\s )})
   = - \stackrel{\cdot \cdot}{\gamma} (t) - \frac{1}{4\pi^2} \frac{\partial V}{\partial \gamma}
   (t,\gamma {\s (}t{\s )}) \\
  &=& -\stackrel{\cdot \cdot}{\gamma} (t) - \frac{1}{2\pi} \sin 2\pi (\gamma {\s (}t{\s )}-kt-q_0) \; . \nonumber
\end{eqnarray}
We observe that
\begin{eqnarray*}
  \gamma^- (t) &=& \gamma_{k,q_0} (t) =kt+q_0 \\
  \gamma^+ (t) &=& \gamma_{k,q_0} (t+\frac{1}{2k} ) =kt+q_0+\frac{1}{2}
\end{eqnarray*}
are solutions of (\ref{eq:perM}) corresponding to the minimum and maximum of $h$.
Equation (\ref{eq:perM}) may be interpreted as describing a mathematical pendulum
with gravity, where the observer rotates with angular velocity $k$. The obvious 
equilibrium states for $k=0$ 'pendulum up' (unstable) and 'pendulum down' (stable)
correspond to $\gamma^-$ and $\gamma^+$. This also holds for general $k$. For the rotating 
observer however these equilibrium states now appear as rotations.
A short calculation shows that ${\cal I}_V(\gamma^- )=2\pi^2 k^2 +1$, ${\cal I}_V(\gamma^+ )=2\pi^2 k^2 -1$
and the Morse index of $\gamma^-$ resp. $\gamma^+$ (regarded as critical points of ${\cal I}_V$)
is $1$ resp. $0$: The eigenvalues of the perturbed Jacobi operator acting on
$C^\infty({\mathbb R}/{\mathbb Z},{\mathbb R})$
\[
  L_V^\mp \xi = - \stackrel{\cdot \cdot}{\xi} - \frac{1}{4\pi^2} \frac{\partial^2 V}{\partial u^2}
   (t,u=\gamma^\mp) \xi =  - \stackrel{\cdot \cdot}{\xi} \mp \xi = \lambda_l^\mp \xi
\]
are given by
\[
  \lambda_l^\mp =4\pi^2 l^2 \mp 1\; ,\; l\in {\mathbb N}_0 \; .
\]
Hence $L_V^+$ has only strictly positive eigenvalues, whereas $L_V^-$ has
exactly one negative eigenvalue of multiplicity one, namely $\lambda_0^- =-1$.
Proposition 2.2 in \cite{CFHW} states that $\gamma^\mp$ are the only solutions
of (\ref{eq:perM}) in $\Lambda_kS^1$, i.e. $Crit\: {\cal I}_V=\{ \gamma^- ,\gamma^+ \}$
is discrete. Therefore the Bott-type Morse complex reduces to the Morse-Witten
complex (cf. \cite{Sch},cite{W93},\cite{W95})
\begin{eqnarray*}
  CM_{i}(\Lambda_k S^1,{\cal I}_V,g;{\mathbb Z}_{2}) & = &
  \begin{cases}                                 
    {\mathbb Z}_2 <\gamma^+ > & ,\; i=0 \: ,\\
    {\mathbb Z}_2 <\gamma^- > & ,\; i=1 \: ,\\
    0 & ,\text{else} \: .                       
  \end{cases}
\end{eqnarray*}   
The only a priori nontrivial matrix element has the coefficient $n_2(\gamma^- ,\gamma^+ )$,
which is defined to be the number of connecting orbits modulo $2$, i.e. solutions
$w\in C^\infty({\mathbb R}\times {\mathbb R}/{\mathbb Z},{\mathbb R})$ of
\begin{equation} \label{eq:gradE}
  \partial_s w-\partial_t \partial_t w-\nabla V(t,w)=0\; ,\; w(s,\cdot )
  \stackrel{s\rightarrow \mp \infty}{\longrightarrow} \gamma^\mp (\cdot ) \; .
\end{equation}
Here we identify two solutions $w,\tilde{w}$, if there exists $\tau \in {\mathbb R}$ such
that $w(s+\tau ,t)=\tilde{w} (s,t)$ for all $s,t$.
The {\sl Ansatz} $w(s,t)=kt+q_0+\chi (s)$, where $\chi \in C^\infty({\mathbb R},{\mathbb R}/{\mathbb Z})$,
leads to the ODE
\begin{equation} \label{eq:gradchi}
  \chi^\prime (s)={\s \frac{1}{2\pi}} \sin 2\pi \chi (s) \; ,
\end{equation}
which has stationary solutions for the initial values
\begin{eqnarray*}
  \chi_0 =\chi (s_0=0)={\s \frac{1}{2}}\; &:& \; \chi (s) \equiv {\s \frac{1}{2}} \; , \\
  \chi_0 =\chi (s_0=0)=0 \; &:& \; \chi (s) \equiv 0 \; . 
\end{eqnarray*}
Choosing initial values $\chi_0 \in (0,{\s \frac{1}{2}})$ resp. $\chi_0 \in ({\s \frac{1}{2}},1)$
$\chi (s)$ behaves as follows
\[
  \chi (s) \stackrel{s\rightarrow -\infty}{\longrightarrow} 0 \; , \; 
  \chi (s) \stackrel{s\rightarrow +\infty}{\longrightarrow} {\s \frac{1}{2}} \; ,
\]
resp.
\[
  \chi (s) \stackrel{s\rightarrow -\infty}{\longrightarrow} 1 \; , \; 
  \chi (s) \stackrel{s\rightarrow +\infty}{\longrightarrow} {\s \frac{1}{2}} \; ,
\]
showing that our {\sl Ansatz} yields two connecting orbits between $\gamma^-$ and $\gamma^+$.
As $Ind_{{\cal I}_V}\: (\gamma^- )=1$ there are no others. 
As a consequence $n_2(\gamma^- ,\gamma^+ )=0$
and the generators of the chain complex coincide with the ones of the corresponding 
homology groups.
Note that equation (\ref{eq:gradchi}) coincides modulo a factor with the gradient flow equation
of the Morse function $h$ on ${\mathbb R}/{\mathbb Z}$.
For general $n$ the (nondegenerate) Morse homology groups coincide with (\ref{eq:Morse-hom}).

Now we treat the case of Floer homology of $T^*S^1$ by perturbing the free Hamiltonian by the
time-dependent potential $V$ as above:
\[
  H(t,u,v)={\s \frac{1}{2}} ({\s \frac{v}{2\pi}})^2 -\cos 2\pi (u-kt-u^0) \;
\]
We fix a solution $x_{k,x^0}(t)=(u_{k,u^0}{\s (}t{\s )},v^0)=(kt+u^0,4\pi^2 k)$
of the unperturbed problem (\ref{6}). Our equation of interest now reads
\begin{eqnarray} \label{eq:ham}
   \begin{pmatrix} \dot u(t) \\ \dot v(t) \end{pmatrix}
  &=&\dot x(t)
   =X_H(x{\s (}t{\s )}) \\
  &=& \begin{pmatrix} v(t)/(2\pi )^2 \\ -\frac{\partial V}{\partial u} (t,u{\s (}t{\s )}) 
    \end{pmatrix}
   =\begin{pmatrix} v(t)/(2\pi )^2 \\ -2\pi \sin 2\pi (u{\s (}t{\s )}-kt-u^0) 
    \end{pmatrix} \; . \nonumber
\end{eqnarray}
We have two solutions
\[
  x^-(t)=x_{k,x^0}(t)=\begin{pmatrix} kt+u^0 \\ 4\pi^2 k \end{pmatrix} \; , \;
  x^+(t)=x_{k,x^0}(t+\frac{1}{2k} )=\begin{pmatrix} kt+u^0+\frac{1}{2} \\ 4\pi^2 k 
  \end{pmatrix} \; , \;
\]
which according to \cite{CFHW} Proposition 2.2 are the only ones. Linearizing (\ref{eq:ham})
at $x^\mp$ yields
\begin{equation} \label{eq:linham}
  \dot \xi (t)=-J\nabla_{\xi (t)} \nabla H(x^\mp ) =-J S^\mp \xi (t) \; ,
\end{equation}
where
\[
  S^\mp = \begin{pmatrix} \frac{1}{4\pi^2} \frac{\partial^2 H}{\partial u^2} (t,x^\mp )
   & \frac{1}{4\pi^2} \frac{\partial^2 H}{\partial u\partial v} (t,x^\mp )
  \\  
   4\pi^2 \frac{\partial^2 H}{\partial v \partial u} (t,x^\mp )
   & 4\pi^2 \frac{\partial^2 H}{\partial v^2} (t,x^\mp ) \end{pmatrix}
  =\begin{pmatrix} \pm 1 & 0 \\ 0 & 1 \end{pmatrix}
\]
and
\[
  J=\begin{pmatrix}  0 & -g^{11} \\  g_{11} & 0 \end{pmatrix}
  =\begin{pmatrix}  0 & \frac{-1}{4\pi^2} \\  4\pi^2 & 0 \end{pmatrix} \; .
\]
The flow given by (\ref{eq:linham}) is a path 
$\Phi^\mp :[0,1]\rightarrow Sp(2,{\mathbb R})$ starting
at the identity: $\Phi^\mp (t)=e^{-tJS^\mp}$. By \cite{SZ} Theorem 3.3 (iv) the Conley-Zehnder
index of $\Phi^\mp$ is given by
\[
  \mu_{CZ} (x^i)=\mu_{CZ} (\Phi^i )=\nu (S^i)-n=
  \begin{cases}                                 
    -1 & ,\; i=- \: ,\\
    0 & ,\; i=+ \: ,                       
  \end{cases}
\]
where $\nu^- (S^i)$ denotes the number of negative eigenvalues of $S^i$ and $n=1$.
Therefore
\[
  \mu_{CZ} (x^\mp )=-Ind_{{\cal I}_V}\: (\gamma^\mp ) ;
\]
a result which has been established in \cite{W99} in the general
nondegenerate case.
The construction of the Floer
cochain complex proceeds as above. Note that it is graded by minus the Morse index
and its cohomology has one generator in dimension $-1$ and one in dimension $0$.
For general $n$ the Floer cohomology coincides with (\ref{eq:Floer-Bott}).

\end{appendix}




\begin{thebibliography}{9999}
    \bibitem[AB]{AB} Austin D.M., Braam P.J., {\sl Morse-Bott theory and equivariant
     cohomology} in {\sl The Floer memorial volume}, PM {\bf 133}, Birkh\"auser 1995. 
    \bibitem[CFHW]{CFHW} Cieliebak K., Floer A., Hofer H., Wysocki C.,
     {\sl Applications of symplectic homology II: stability of the action spectrum},
     preprint ETH Z\"urich.
    \bibitem[CZ]{CZ} Conley C., Zehnder E., {\sl Morse type index theory for flows
     and periodic solutions for Hamiltonian equations}, Comm. Pure Appl. Math.
     {\bf XXXVII} (1984), 207-253.
    \bibitem[Fl]{Fl} Floer A., {\sl Symplectic fixed points and holomorphic
     spheres}, Comm. Math. Phys. {\bf 120} (1989), 575-611.
    \bibitem[GL]{GL} Gelfand I.M., Lidskii V.B., {\sl On the structure 
     of the regions of stability of linear canonical systems of differential 
     equations with periodic coefficients},
     Translations A.M.S. ({\bf 2}) 8 (1958), 143-181.
    \bibitem[Jo]{Jo} Jost J., {\sl Riemannian geometry and geometric analysis},
     Universitext, Springer-Verlag 1995.
    \bibitem[Kl]{Kl} Klingenberg W., {\sl Lectures on closed geodesics},
      Grundlehren der mathematischen Wissenschaften {\bf 230}, Springer-Verlag 1978.
    \bibitem[RS]{RS} Robbin J., Salamon D., {\sl The Maslov index for paths},
      Topology {\bf 32} (1993), 827-844.
    \bibitem[RT]{RT} Ruan Y., Tian G., {\sl Bott-type symplectic Floer cohomology 
      and its multiplication structures}, Math. Res. Letters {\bf 2} (1995), 203--219.
    \bibitem[Sch]{Sch} Schwarz M., {\sl Morse homology}, PM {\bf 111}, Birkh\"auser 1993.
    \bibitem[SW]{SW} Salamon D., Weber J., {\sl $J$-holomorphic curves in cotangent
     bundles and Morse theory on the loop space}, in preparation.
    \bibitem[SZ]{SZ} Salamon D., Zehnder E., {\sl Morse theory for periodic solutions
     of Hamiltonian systems and the Maslov index}, Comm. Pure Appl. Math. {\bf XLV}
     (1992), 1303-1360.
    \bibitem[Vi]{Vi} Viterbo C., {\sl Functors and computations in Floer homology with
     applications -- part II}, preprint October 1996.
    \bibitem[W93]{W93} Weber J., {\sl Der Morse-Witten Komplex}, Diplomarbeit am FB Mathematik
     der TU Berlin, Februar 1993.
    \bibitem[W95]{W95} Weber J., {\sl Morse-Ungleichungen, Supersymmetrie und quasiklassischer
     Limes}, Diplomarbeit am FB Physik der TU Berlin, Februar 1995.
    \bibitem[W98]{W98} Weber J., {\sl Topology of $Sp(2,{\mathbb R})$ and the Conley-
     Zehnder index}, preprint University of Warwick 51/1998.
    \bibitem[W99]{W99} Weber J., {\sl $J$-holomorphic curves in cotangent bundles
     and the heat flow}, PhD-thesis TU Berlin 1999.
  \end{thebibliography}
\end{document}